\pdfoutput=1
\documentclass[preprint,12pt]{elsarticle}

\usepackage{amssymb}
\usepackage{amsmath}
\usepackage{diagbox}  
\usepackage{bm}%
\usepackage[linesnumbered,ruled,vlined]{algorithm2e}
\usepackage{subcaption}
\usepackage{booktabs}%

\newtheorem{remark}{Remark}%
\newtheorem{example}{Example}%
\newcommand{\argmin}{\mathop{\text{argmin}}}

\journal{}

\begin{document}

\begin{frontmatter}



\title{Alternately-optimized SNN method for acoustic scattering problem in unbounded domain}

\author[a]{Haoming Song}
\ead{songhaoming314@163.com}
\author[b]{Zhiqiang Sheng}
\ead{sheng_zhiqiang@iapcm.ac.cn}
\author[c]{Dong Wang}
\ead{wangdong@cuhk.edu.cn}
\author[a]{Junliang Lv\corref{cor1}}
\ead{lvjl@jlu.edu.cn}
\cortext[cor1]{Corresponding author: Junliang Lv, lvjl@jlu.edu.cn}

\affiliation[a]{
	organization={School of Mathematics, Jilin University},
	city={Changchun},
	postcode={130012},
	country={China}
}
\affiliation[b]{
	organization={National Key Laboratory of Computational Physics, Institute of Applied Physics and Computational Mathematics},
	city={Beijing},
	postcode={100088},
	country={China}
}
\affiliation[c]{
	organization={School of Science and Engineering, The Chinese University of Hong Kong},
	 city={Shengzhen},
	postcode={518172},
	country={China}
}

\begin{abstract}
In this paper, we propose a novel machine learning-based method to solve the acoustic scattering problem in unbounded domain. 
We first employ the Dirichlet-to-Neumann (DtN) operator to truncate the  physically unbounded domain into a computable bounded domain. 
This transformation reduces the original scattering problem in the unbounded domain to a boundary value problem within the bounded domain.
To solve this boundary value problem, we design a neural network with a subspace layer, where each neuron in this layer represents a basis function. 
Consequently, the approximate solution can be expressed by a linear combination of these basis functions. 
Furthermore, we introduce an innovative alternating optimization technique which alternately updates the basis functions and their linear combination coefficients respectively by training and least squares methods.
In our method, we set the coefficients of basis functions to 1 and use a new loss function each time train the subspace.
These innovations ensure that the subspace formed by these basis functions is truly optimized.
We refer to this method as the alternately-optimized subspace method based on  neural networks (AO-SNN).
Extensive numerical experiments demonstrate that our new method can significantly reduce the relative $l^2$ error to $10^{-7}$ or lower, outperforming existing machine learning-based methods to the best of our knowledge.
\end{abstract}



\end{frontmatter}



\section{Introduction}
The rapid development of machine learning has attracted researchers to apply them in solving various scientific and engineering problems. 
Due to the powerful nonlinear approximation capabilities \cite{Approximator1, Approximator2} of neural networks, machine learning-based methods have been found to be effective in approximating the solutions of many partial differential equations. 
These machine learning-based methods include the physics-informed neural networks (PINN) \cite{PINN_2019}, the deep Ritz method (DRM) \cite{DRM_2018}, the weak adversarial networks (WAN) \cite{WAN_2020} and the deep finite volume method (DFVM) \cite{DFVM_2024}. The basic idea of these methods is to employ a properly designed neural network to approximate the solution of the PDEs, and then optimize all parameters of the neural network simultaneously through the gradient descent methods, such as the stochastic gradient descent (SGD) and the Adam method. One of the differences among these methods is that they have different loss functions. The PINN adopts a loss function defined by the strong form of the PDEs at collocation points. The loss function of DRM is designed from the principle of minimal potential energy. The loss function used by the WAN is constructed using the weak form of the equation, where two neural networks are introduced to parameterize the weak solution and the test function. The DFVM uses the local conservation law of the flux to form its loss function. 

Unlike the previously mentioned methods that directly represent the numerical solution using the output of the neural network, other methods first use the neural network to construct a set of basis functions and then use a linear combination of them to represent the numerical solution. 
In \cite{SNN_arXiv}, the authors propose the subspace method based on the neural networks (SNN) that trains a neural network to form a subspace with effective representation capabilities for solution of the PDEs and then seeks an approximate solution in this subspace by the least squares method. In \cite{RobustTrain_2020}, the authors proposed the least squares/gradient descent (LSGD) method, which alternately used the gradient descent method to optimize the basis functions and used the least squares methods to optimize the coefficients.
Both SNN and LSGD can achieve higher accuracy than the PINN. To reduce the computational cost, the authors of \cite{RFM_2022} proposed the random feature method (RFM) that used a shallow neural network with randomly generated parameters to form a set of bases. Recently, the transferable neural networks (TransNet) was proposed in \cite{TransNet_2024} to enhance the representation capabilities of the randomly generated basis functions by uniformly distributing the partition hyperplane of the basis functions and tuning the shape parameters of the basis functions. The TransNet can use one pre-trained neural network to solve several PDE problems.

Scattering problems are very important in science and engineering. 
Numerically solving scattering problems has many challenges, such as complex-valued solutions, unbounded physical regions, and oscillations in high-frequency problems. Some special treatments are often required in numerically solving scattering problems, even if machine learning-based methods are used.
For the acoustic scattering problem in bounded domain,
the authors of \cite{PWNN_2022} designed the plane wave activation based neural network (PWNN) to efficiently solve the Helmholtz equation. 
In \cite{RBDNN_2022}, the ray-based deep neural network (RBDNN) was presented by using a combination of deep neural network and the plane wave function to approximate the solution of the acoustic wave scattering problem. 
Both PWNN and RBDNN can solve high-frequency problems efficiently.
For the acoustic scattering problem in unbounded domain, neither traditional nor machine learning-based methods can directly solve this type of problem. 
To transform the problems in unbounded physical domain into the problems in bounded computational domain, researchers have proposed many domain truncation techniques such as the absorbing boundary condition (ABC) \cite{ABC_1977}, the perfectly matched layer method (PML) \cite{PML_1994,CW2003}, the transparent boundary condition (TBC) \cite{TBC_2021,TBC_2024,TBC_2025} and the perfectly absorbing layer method (PAL) \cite{PAL_2019}.
Recently, researchers have started to use machine learning-based methods combined with domain truncation techniques to solve problems in unbounded domain.
In \cite{RBDNN_2023}, the authors proposed the RBDNN-PML method, which first truncates the physical domain by using the PML method and then solves the PML problem in the bounded domain by using the RBDNN method.
The authors of \cite{cPINN_2024} combined the PINN and the first-order ABC methods to design the coupling physics-informed neural networks (cPINN) method to solve the acoustic wave scattering problems in unbounded domains and obtained some convergence results.
However, the methods mentioned are difficult to obtain high accuracy.

In the present paper, we attempt to construct a machine learning-based method with high accuracy to solve the acoustic wave scattering problems in unbounded domains. 
Due to the TBC method's ability to precisely transform problems from unbounded domains to bounded ones, we use it to truncate the unbounded domain, thereby obtaining a boundary value problem.
However, if we use the standard SNN method directly to solve this boundary value problem, we find that the numerical solution only achieves an accuracy of about $10^{-3}$. 
The reason for this degradation in accuracy may be the inability of SNN methods to obtain a good enough subspace in solving complex scattering problems with only one training.
To obtain a good enough subspace, an obvious strategy may be to alternately update the subspace and the coefficients of linear combinations. 
The method derived along this line is actually close to the LSGD \cite{RobustTrain_2020}. 
However, this method still fails to obtain high accuracy in solving the boundary value problem transformed from the scattering problem in unbounded domain. 
In fact, numerical experiments show that the magnitude of the coefficients of the linear combinations obtained by this method always varies considerably.
This phenomenon indicates that the role played by these basis functions in representing the approximate solution varies significantly.

In order to make the subspace truly optimized, we propose a new alternating optimization method called the alternatively-optimized subspace method based on neural networks (AO-SNN).
There are two major innovations in the AO-SNN.
First, we set the linear combination coefficients of the basis functions to 1 during each training. 
This strategy can be interpreted as a form of normalization, which ensures that the roles of the resulting basis functions in representing the approximate solution are comparable.
Moreover, we design a new loss function to measure the difference between the currently trained neural network and the most recent approximate solution.
After obtaining new basis functions in this way, we then use the least squares method to obtain a new approximate solution.
Performing the above procedure several times alternately, we can obtain a numerical solution with high accuracy, about $10^{-7}$.

The rest of this paper is organized as follows. 
In Section 2, we introduce the model problem of the acoustic scattering in unbounded domain and use the DtN operator to reformulate it to a boundary value problem in bounded domain.  
In Section 3, the AO-SNN method is proposed to solve this boundary value problem. In Section 4, we give some numerical examples and compare the AO-SNN method with other methods in accuracy. Then, some discussions about the AO-SNN method are given in Section 5. In the last section, we give some conclusions.

\section{Model problem}
\subsection{Acoustic scattering in unbounded domain}

\begin{figure}[htp!]
	\begin{center}
		\includegraphics[scale=0.25]{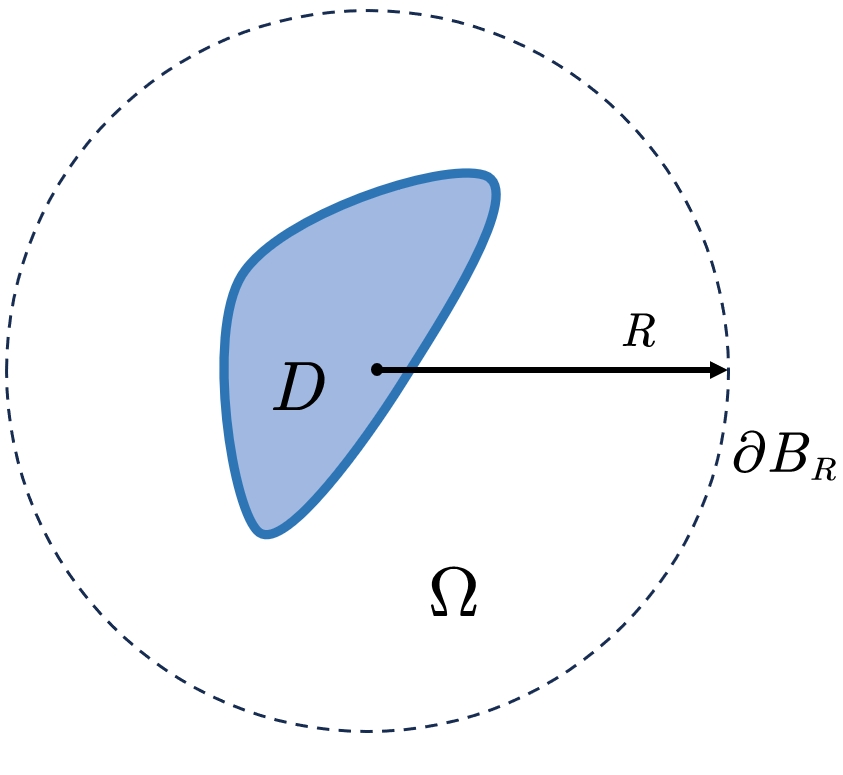}
	\end{center}
	\caption{Geometry of the obstacle scattering problem}\label{AS}
\end{figure}

We consider the two-dimensional scattering of time-harmonic acoustic waves by an obstacle $D$ with 
Lipschitz boundary $\partial D$; see Figure \ref{AS}. Let $B_R = \{ \bm{x}\in\mathbb{R}^2:|\bm{x}|<R\}$ be a ball with radius $R$, where $R$ is large enough such that $\Bar{D}\subset B_{R}$. 


The process of acoustic wave 
propagation can be governed by the Helmholtz equation
\begin{equation}
	\label{Heq}\Delta u + \kappa^2u = 0, ~~~~{\rm in}~\mathbb{R}^2\setminus \bar{D},
\end{equation}
where $\kappa>0$ is called the wavenumber.  Depending on the material of the obstacle D,
the boundary conditions on $\partial D$ can be divided into three types as follows.
\begin{itemize}
	\item Sound soft (Dirichlet boundary condition)
	\begin{equation}\label{ssbc}
		u=g(\bm{x}),~~~~{\rm on}~\partial D.
	\end{equation} 
	\item Sound hard (Neumann boundary condition)
	\begin{equation}\label{shbc}
		\frac{\partial u}{\partial \bm{n}}=g(\bm{x}),~~~~{\rm on}~\partial D,
	\end{equation} 
	where $\frac{\partial}{\partial \bm{n}}$ means the normal derivative on $\partial D$.
	\item Impedance boundary condition (Robin boundary condition)
	\begin{equation}\label{imbc}
		\frac{\partial u}{\partial \bm{n}} + \mathrm{i}\lambda u = g(\bm{x}),~~~~{\rm on}~\partial D,
	\end{equation} 
	where $\lambda$ is the impedance coefficient.
\end{itemize}
In addition, to ensure the uniqueness of the solution, the scattered wave field $u$ is required to satisfy the Sommerfeld radiation condition at infinity
\begin{equation}
	\label{SRC}\lim_{r\to\infty}\sqrt{r}(\frac{\partial u}{\partial r}-\mathrm{i}\kappa u)=0, ~~~~r=|\bm{x}|,
\end{equation}
where $\frac{\partial}{\partial r}$ means the partial derivative with respect to $r$.

\subsection{The transparent boundary condition}

Since the physical domain of the scattering problem is unbounded, we need to transform this problem into a boundary value problem in bounded domain when we use classical numerical methods, such as the finite element method \cite{JLLWWZ2022,JLLZ2017_ESAIM,JLLZ2017_JSC} and the finite difference method. 
Often used truncation techniques include the absorbing boundary condition \cite{ABC_1977}, the transparent boundary condition \cite{TBC_2021}, the perfectly matched layer method \cite{CW2003}, and the perfectly absorbing layer method \cite{PAL_2019}.
Since the TBC method can accurately truncate the scattering problem in unbounded domain into an boundary value problem in bounded domain, we choose the TBC method in this paper.

In the exterior domain $\mathbb{R}^2\setminus\bar{B}_{R}$, the solution $u$ of 
equation (\ref{Heq}) can be expressed as a Fourier expansion in polar coordinates $(r,\theta)$
\begin{align}
	u(r,\theta) &= \sum_{n\in\mathbb{Z}}\frac{H^{(1)}_n(\kappa r)}{H^{(1)}_n(\kappa R)}\hat{u}_ne^{\mathrm{i}n\theta},\\
	\label{uhat}\hat{u}_n &= \frac{1}{2\pi}\int_0^{2\pi}u(R,\theta)e^{-\mathrm{i}n\theta}\mathrm{d}\theta,~~~~\forall r>R,
\end{align}
where $H^{(1)}_n(\cdot)$ is the Hankel function of the first kind with order $n$.

Given $u$ belonging to the trace space $H^{1/2}(\partial B_R)$, 
we have the DtN operator $\mathcal{F}:H^{1/2}(\partial B_R)\to H^{-1/2}(\partial B_R)$, defined by
\begin{equation}
	\label{DtN}\mathcal{F}[u](R,\theta) = \frac{1}{R}\sum_{n\in\mathbb{Z}}h_n(\kappa R)\hat{u}_ne^{\mathrm{i}n\theta},
\end{equation}
where
\begin{equation}
	\label{hn}h_n(z) = z \frac{H^{(1)'}_n(z)}{H^{(1)}_n(z)}.
\end{equation}
Combining (\ref{uhat})-(\ref{hn}), we get the transparent boundary condition
\begin{equation}\label{TBC}
	\mathcal{F}[u] = \frac{\partial u}{\partial r}, ~~~~{\rm on}~\partial B_R.
\end{equation}

\subsection{Boundary value problem with TBC}
By using the transparent boundary condition, we can reformulate acoustic scattering problem (\ref{Heq})-(\ref{SRC}) into a boundary value problem in bounded domain $\Omega=B_{R}\setminus \bar{D}$
\begin{equation}\label{ASwithTBC}
	\begin{cases}
		\Delta u + \kappa^2u = 0, ~~~~&{\rm in}~\Omega,\\
		\mathcal{B} u = g(\bm{x}), ~~~~&{\rm on}~\partial D,\\
		\mathcal{F}[u] = \frac{\partial u}{\partial r}, ~~~~&{\rm on}~\partial B_R,
	\end{cases}
\end{equation}
where $\mathcal{B}$ is the corresponding partial differential operator in equation (\ref{ssbc}), (\ref{shbc}), or (\ref{imbc}).

In equation (\ref{DtN}), the DtN operator in the TBC is defined by an infinite series. In practical computation, it is necessary to truncate the DtN operator by a sufficiently large $N$ and yield the truncated DtN operator
\begin{equation}\label{TBCN}
	\mathcal{F}_N[u](R,\theta) = \frac{1}{R}\sum^{N}_{n=-N}h_n(\kappa R)\hat{u}_ne^{\mathrm{i}n\theta}.
\end{equation}
In addition, $\hat{u}_n$ can be approximated by applying the numerical integration formula to the integral in equation (\ref{uhat}). 

By using the operator $\mathcal{F}_N$, we get the approximate problem to problem (\ref{ASwithTBC})
\begin{equation}\label{ASwithTTBC}
	\begin{cases}
		\Delta u + \kappa^2u = 0, ~~~~&{\rm in}~\Omega,\\
		\mathcal{B} u = g(\bm{x}), ~~~~&{\rm on}~\partial D,\\
		\mathcal{F}_{N}[u] = \frac{\partial u}{\partial r}, ~~~~&{\rm on}~\partial B_R.
	\end{cases}
\end{equation}

\section{The alternately-optimized SNN method}
In this section, we develop the alternately-optimized subspace method based on neural networks (AO-SNN) to solve the approximate problem (\ref{ASwithTTBC}).
Specifically, in Subsection 3.1, we give an introduction of the SNN. 
In Subsection 3.2, we propose the AO-SNN methods to solve general PDEs.
We then give an explain in Subsection 3.3 how the AO-SNN method optimizes the subspace and obtains an approximate solution.
Finally, We give the specific implementations of AO-SNN for solving the approximate problem (\ref{ASwithTTBC}) in Subsection 3.4.

\subsection{The subspace method based on neural network}

In this subsection, we use the two-dimensional linear differential equations
\begin{equation}\label{target_equation}
	\begin{cases}
		\mathcal{A}u=f(u),&\rm{in}~\Omega,\\
		\mathcal{B}u=g, &\rm{on}~\partial\Omega
	\end{cases}
\end{equation}
as an example to describe the stander SNN method. Here, the differential operator $\mathcal{A}$ is a linear operator (e.g. Laplace operator), and the operator $\mathcal{B}$ represents 
the identity operator (corresponding to Dirichlet boundary condition),
the normal derivative operator (corresponding to Neumann boundary condition), or mixed derivative operator (corresponding to Robin boundary condition).

The key idea of the SNN method is to use neural networks to construct a subspace
$$
\mathcal{S}={\rm span}\{\phi_j(\bm{x};\bm{\theta})\}_{j=1}^M, 
$$
where $\{\phi_j(\bm{x};\bm{\theta})\}_{j=1}^M$ is a class of basis functions defined by a neural network with parameters $\bm{\theta}$. Once the subspace $\mathcal{S}$ has been adequately trained, the numerical solution can be expressed as a linear combination $\sum_{j=1}^{M}\omega_j\phi_j(\bm{x};\bm{\theta})$, and determined by solving a least square problem.

\begin{figure}[htp!]
	\begin{center}
		\includegraphics[scale=0.25]{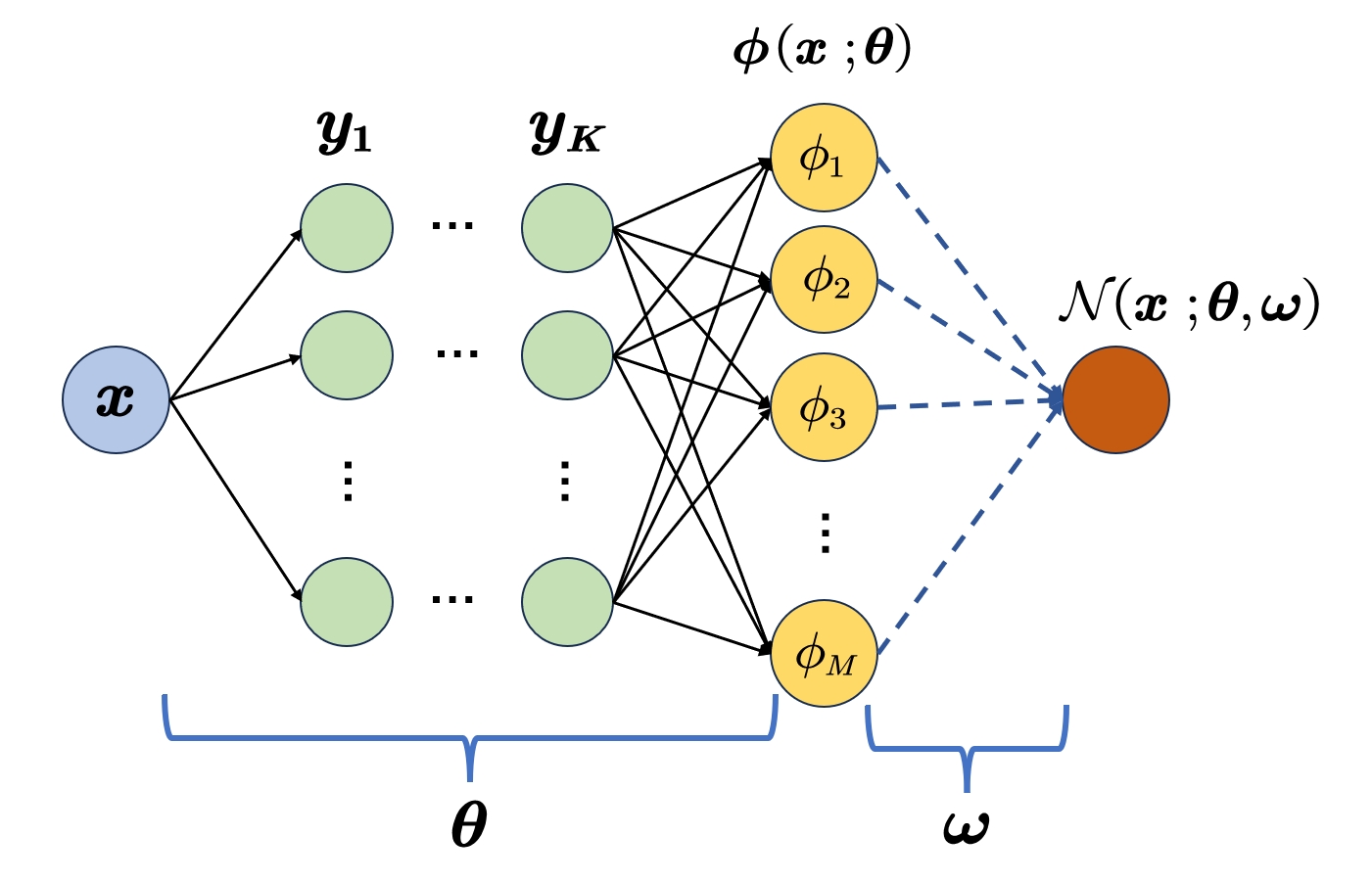}
	\end{center}
	\caption{The architecture of SNN network}
	\label{SNN_network}
\end{figure}

A straightforward way to implement this idea may be to define a neural network consisting of an input layer, several hidden layers, a subspace layer, and an output layer without bias and activation functions, as in Figure \ref{SNN_network}. 
In this architecture, each neuron in the subspace layer can be understood as a function $\phi_j(\bm{x};\bm{\theta})$. These basis functions can span a finite-dimensional subspace $\mathcal{S}$.

Let $K$ represent the number of hidden layers, and $n_j, j=1,\cdots,K$ represent the number of neurons in $j$th hidden layer. Let $M$ be the number of neurons in subspace layer, and $\bm{\phi}=(\phi_1, \phi_2,\cdots, \phi_M)^T$ represent the basis function vector. Let $\bm{\omega}=(\omega_1, \omega_2, \cdots, \omega_M)^T$ represent the weights of last layer (also the coefficients of linear combinations $\mathcal{N}(\bm{x}; \bm{\theta}, \bm{\omega}) = \sum_{j=1}^{M}\omega_j\phi_j$). 
In this way, the forward propagation process of the SNN can be expressed as
\begin{equation}\label{forward_SNN}
	\begin{cases}
		\bm{y}_0 = \bm{x},&\\
		\bm{y}_k = \sigma_k(\bm{W}_k\bm{y}_{k-1} + \bm{b}_k),&k=1,2,\cdots, K+1,\\
		\bm{\phi}(\bm{x}; \bm{\theta}) = \bm{y}_{K+1},&\\
		\mathcal{N}(\bm{x}; \bm{\theta}, \bm{\omega}) = \sum\limits_{j=1}^{M}\omega_j\phi_j(\bm{x}; \bm{\theta}),&
	\end{cases} 
\end{equation}
where $\bm{W}_k\in \mathbb{R}^{n_k\times n_{k-1}}$ and $\bm{b}_k\in \mathbb{R}^{n_k}$ denote the weights and biases of the hidden layers respectively and $\sigma_k(\bm{x})$ is the activation function. Let $\bm{x}$ represent the input value, $n_0$ represent the dimension of $\bm{x}$, and $\bm{\theta} = \{\bm{W}_1, \cdots, \bm{W}_{K+1}, \bm{b}_1, \cdots, \bm{b}_{K+1}\}$ represent all parameters in the hidden layers. Thus, the neural network subspace $\mathcal{S}$ can be defined as 
\begin{equation*}
	\mathcal{S}(\bm{x};\bm{\theta})={\rm span}\{\phi_1, \cdots, \phi_M\}, 
\end{equation*}
and the neural network solution can be represented by $$\mathcal{N}(\bm{x}; \bm{\theta}, \bm{\omega})= \sum_{j=1}^{M}\omega_j\phi_j(\bm{x}; \bm{\theta}).$$ 

There are two main steps in the implementation of the SNN method: training the subspace and solving the least squares problem. The parameters $\bm{\theta}$ and $\bm{\omega}$ are determined in these two steps, respectively.
To train the subspace, we adopt the loss function \begin{equation}\label{general_lossfunction}
	\begin{aligned}
		\mathcal{L}(\mathcal{N}(\bm{x}; \bm{\theta}, \bm{\omega})) &= \frac{1}{N}\sum^{N}_{i=1}\left(\mathcal{A}\big(\mathcal{N}(x_i; \bm{\theta}, \bm{\omega})\big)-f(x_i)\right)^2 \\
		&\quad+ \lambda \frac{1}{\bar{N}}\sum^{\bar{N}}_{i=1}\left(\mathcal{B}\big(\mathcal{N}(\bar{x}_i; \bm{\theta}, \bm{\omega})\big)-g_i(\bar{x}_i)\right)^2,
	\end{aligned}
\end{equation}
where $x_i$ and $\bar{x}_i$ respectively represent the collocation points located in $\Omega$ and on $\partial\Omega$, and $N$ and $\bar{N}$ respectively represent the total number of these two kinds of collocation points.  
This loss function is based on the strong form of partial differential equations with a boundary penalty term, and often used in the PINN \cite{PINN_2019}. 
Once the subspace is sufficiently trained, the parameter $\bm{\theta}$ and the basis functions $\phi_j(\bm{x};\bm{\theta})$ are obtained. Since the differential operators $\mathcal{A}$ and $\mathcal{B}$ are both linear, the optimal parameter $\bm{\omega}$ can be determined by solving the least square problem 
\begin{equation}\label{general_lstsq}
	\argmin_{\bm{\omega}}	\Big(||\sum^M_{j=1}\omega_j\mathcal{A}(\phi(\bm{x};\bm{\theta})) - f(\bm{x})||_{l^2({X})}^2 + ||\sum^M_{j=1}\omega_j\mathcal{B}(\phi(\bm{\bar{x}};\bm{\theta})) - g(\bm{\bar{x}})||_{l^2(\bar{X})}^2\Big),
\end{equation}
where $X$ and $\bar{X}$ represent the sets of collocation points that are located 
in $\Omega$ and on $\partial\Omega$, respectively. 
Given that we have obtained the parameters $\bm{\theta}$ and $\bm{\omega}$, the numerical solution of the SNN method can be expressed as $\mathcal{N}(\bm{x}; \bm{\theta}, \bm{\omega})$.
The implementation of the SNN method is showed in Algorithm 1.

\begin{algorithm}[htp]
	\caption{Implementation of the SNN}  
	Initialize a neural network $\mathcal{N}(\bm{x}; \bm{\theta}, \bm{\omega})$\;
	Set $\bm{\omega}=\bm{1}$ and train the network with the loss function $\mathcal{L}(u)$  (\ref{general_lossfunction}), until the loss is small enough\;
	Solve the least square problem  (\ref{general_lstsq})\;
	Output the neural network solution $\mathcal{N}(\bm{x}; \bm{\theta}, \bm{\omega})$;\
\end{algorithm} 

It is worth pointing out that the SNN method uses different approaches to optimize the parameters $\bm{\theta}$ and $\bm{\omega}$, while the PINN method optimizes both parameters $\bm{\theta}$ and $\bm{\omega}$ using one approach (e.g. SGD or adam).
It is this strategy employed by the SNN method that allows it to achieve higher accuracy than the PINN.

\subsection{The Alternately-optimized SNN method}\label{AO-SNN}

Although the SNN method is generally more accurate than the standard PINN, it can only achieve an accuracy about $10^{-3}$ in solving boundary value problem (\ref{ASwithTTBC}). 
This loss of accuracy may be due to the inability to obtain good enough basis functions to approximate the solution of complex scattering problem using only one training.
A natural improvement is to alternately optimize parameters $\bm{\theta}$ and $\bm{\omega}$ until the  desired accuracy is reached.
Such an idea has been proposed in \cite{RobustTrain_2020} for solving some standard PDEs.
However, this improvement by itself does not enhance accuracy in solving the boundary value problem (\ref{ASwithTTBC}).
In fact, a noticeable difference in the magnitudes of the parameters $\{\omega_1, \cdots, \omega_M\}$ obtained by using this method can be observed from numerical experiments.
This phenomenon implies that the contribution of each basis function $\phi_j$ in approximating the solution is inconsistent.
Therefore, using only the standard alternating optimization strategy is not really enough to get better basis functions.

\begin{figure}[htp]
	\begin{center}
		\includegraphics[scale=0.21]{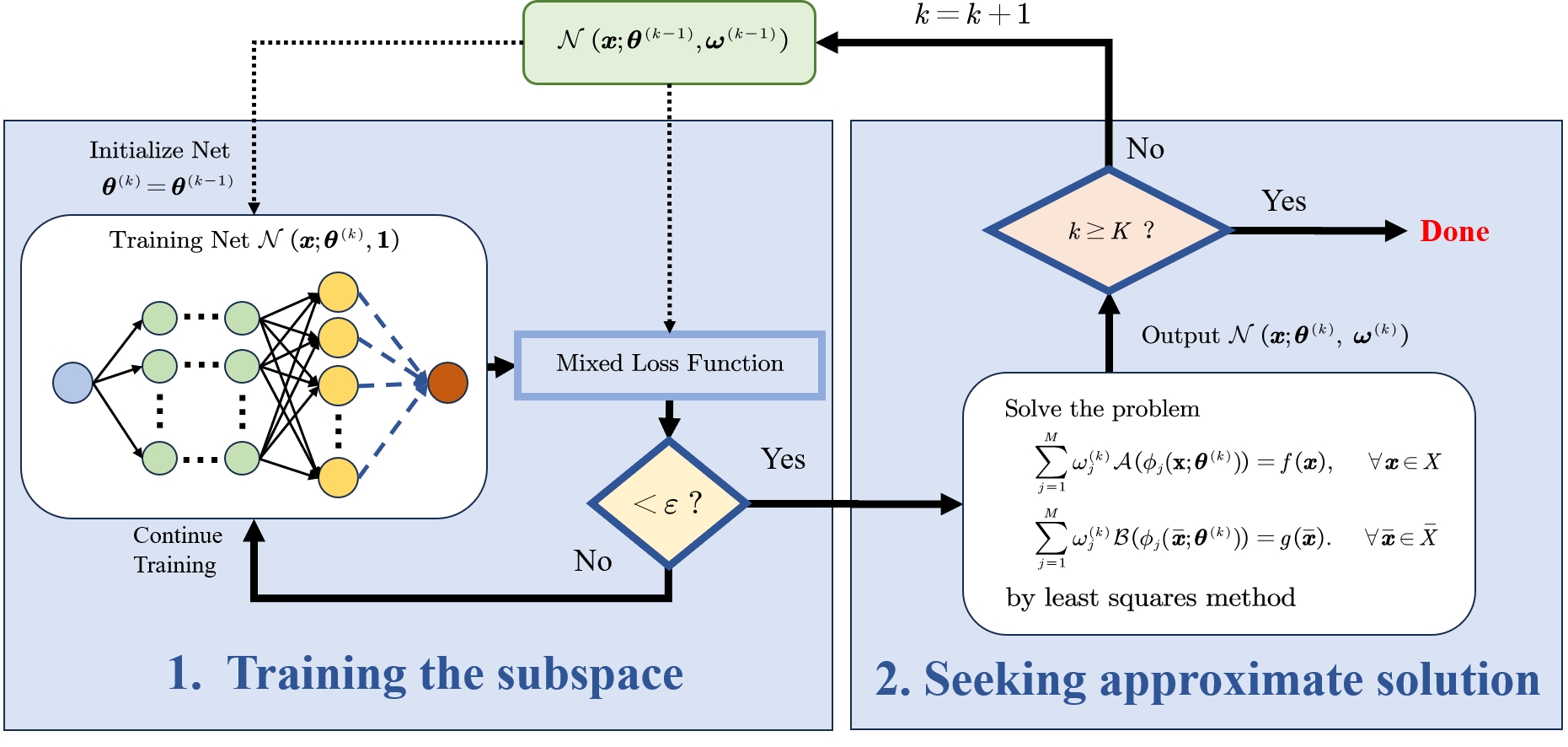}
	\end{center}
	\caption{Flowchart of the AO-SNN}\label{ChartAOSNN}
\end{figure}

In our method, we creatively use a new loss function and set the value of $\bm{\omega}$ during training stages. Specifically, in order to optimize the subspace efficiently in the training stage, i.e., to obtain basis functions that can play nearly the same role in expressing the solution, we set $\bm{\omega}=\bm{1}$ in every training stage.
However, it is worth pointing out that if we still use the loss function (\ref{general_lossfunction}) when training the subspace in this way, we are actually just repeating the SNN method with different initial guesses. 
So we also need to design a new loss function for training. 
In fact, when we apply the SNN method once to the target equations (\ref{target_equation}), we get an approximate solution $\mathcal{N}(\bm{x};\bm{\theta}^{(0)}, \bm{\omega}^{(0)})$.
After setting $\bm{\omega}=\bm{1}$ and initializing  $\bm{\theta}^{(1)}=\bm{\theta}^{(0)}$, we train the neural network to update $\bm{\theta}^{(1)}$, minimizing some loss function.
Thus, we get an intermediate solution $\mathcal{N}(\bm{x}; \bm{\theta}^{(1)}, \bm{1})$.
This loss function should be chosen such that $\mathcal{N}(\bm{x}; \bm{\theta}^{(1)}, \bm{1})$ is an approximate solution to the target equations and is very close to $\mathcal{N}(\bm{x}; \bm{\theta}^{(0)}, \bm{\omega}^{(0)})$.
To this end, we choose the following loss function
\begin{equation}\label{Mixed_loss}
	\begin{aligned}
		\mathcal{L}^{\rm Mixed}\big(\mathcal{N}(\bm{x}; \bm{\theta}^{(k)}, \bm{1})\big) & = \eta_k\mathcal{L}^{\rm Metric}(\mathcal{N}(\bm{x}; \bm{\theta}^{(k)}, \bm{1}), \mathcal{N}(\bm{x}; \bm{\theta}^{(k-1)}, \bm{\omega}^{(k-1)}))\\
		& \quad +\sigma_k\mathcal{L}(\mathcal{N}(\bm{x}; \bm{\theta}^{(k)}, \bm{1})),
	\end{aligned}
\end{equation}
where $\eta_k$ and $\sigma_k$ are hyperparameters, and $\mathcal{L}^{\rm Metric}$ is the metric loss function defined by 
\begin{equation}\label{metric_loss}
	\mathcal{L}^{\rm Metric}\big(\mathcal{N}(\bm{x}; \bm{\theta}, \bm{\omega}), \mathcal{N}(\bm{x}; \bm{\theta}', \bm{\omega}')\big)=\frac{1}{N}\sum^{N}_{i=1}\sum_{|\bm{\alpha}|=0}^{\gamma}\big(\partial^{\bm{\alpha}}\mathcal{N}(x_i; \bm{\theta}, \bm{\omega})-\partial^{\bm{\alpha}}\mathcal{N}(x_i; \bm{\theta}', \bm{\omega}')\big)^2.
\end{equation}
Here, $\gamma$ is a nonnegative integer,  $\bm{\alpha}=(\alpha_1, \alpha_2)$ is a multiple indicator with $\alpha_1$ and $\alpha_2$ being two  nonnegative integers, and $|\bm{\alpha}|=\alpha_1+\alpha_2$. 
Once $\bm{\theta}^{(1)}$ is obtained, we get an intermediate solution $\mathcal{N}(\bm{x}; \bm{\theta}^{(1)}, \bm{1})$ and a new subspace $\mathcal{S}(\bm{x};\bm{\theta}^{(1)})$. 
Using the least squares method, we can then compute $\bm{\omega}^{(1)}$, which produces a new approximate solution $\mathcal{N}(\bm{x}; \bm{\theta}^{(1)}, \bm{\omega}^{(1)})$ in subspace $\mathcal{S}(\bm{x};\bm{\theta}^{(1)})$.
In this way, we have completed an alternating optimization of parameters $\bm{\theta}$ and $\bm{\omega}$. 
Repeat the above process until the final approximate solution is obtained. 
We refer to the this algorithm as the alternately-optimized subspace method based on neural networks (AO-SNN).
The workflow of the AO-SNN method is illustrated in Figure \ref{ChartAOSNN} and the implementation of the AO-SNN method is demonstrated in Algorithm \ref{AlgorithmAOSNN}.

\begin{remark}
	In the mixed loss function $\mathcal{L}^{\rm Mixed}$, the metric loss function $\mathcal{L}^{\rm Metric}$ plays a key role in training the subspace, while the loss function $\mathcal{L}$ enhances the stability of the AO-SNN method. 
	If one chooses $\eta_k=0$ and $\sigma_k=1$, The $k$th AO-SNN iteration happens to be the SNN method with a specific initial guess.
	In this case, high-precision numerical results cannot be achieved.
	In practice, we can usually take $\eta_k=1$ and $\sigma_k=1$.
	However, in numerical experiments, we find that the accuracy and efficiency of the AO-SNN taking $\eta_k=1$ and $\sigma_k=0$ do not differ significantly from the former. 
	This phenomenon further indicates that $\mathcal{L}^{\rm Metric}$ plays a more important role than $\mathcal{L}$ in the mixed loss function.
\end{remark}

\begin{remark}
	When the non-negative integer $\gamma$ in the metric loss function is smaller than the order of the highest derivative in the target PDEs, a larger $\gamma$ leads to higher accuracy in the approximate solution obtained by the AO-SNN method, but also increases the computational cost.
	When the value of $\gamma$ exceeds the order of the highest-order derivative, further increases in $\gamma$ do not lead to greater accuracy.
	In practical calculations, we can choose $\gamma$ to be the order of the highest order derivative in the target PDEs.
\end{remark}

\begin{algorithm}[htp]
	\caption{Implementation of the AO-SNN}  \label{AlgorithmAOSNN}
	Get an initial neural network solution $\mathcal{N}(\bm{x}; \bm{\theta}^{(0)}, \bm{\omega}^{(0)})$ by SNN with loss function (\ref{general_lossfunction})\;  
	
	\For{$k=1,2,\cdots,K$}   
	{  
		Train the network with the mixed loss function (\ref{Mixed_loss}) to get $\bm{\theta}^{(k)}$\;
		Solve the least square problem (\ref{general_lstsq}) to get $\bm{\omega}^{(k)}$\;
	}  
	Output the final approximate solution $\mathcal{N}(\bm{x}; \bm{\theta}^{(K)}, \bm{\omega}^{(K)})$.  
\end{algorithm} 

\subsection{An explanation of the AO-SNN}
\begin{figure}[htp]
	\begin{center}
		\includegraphics[scale=0.23]{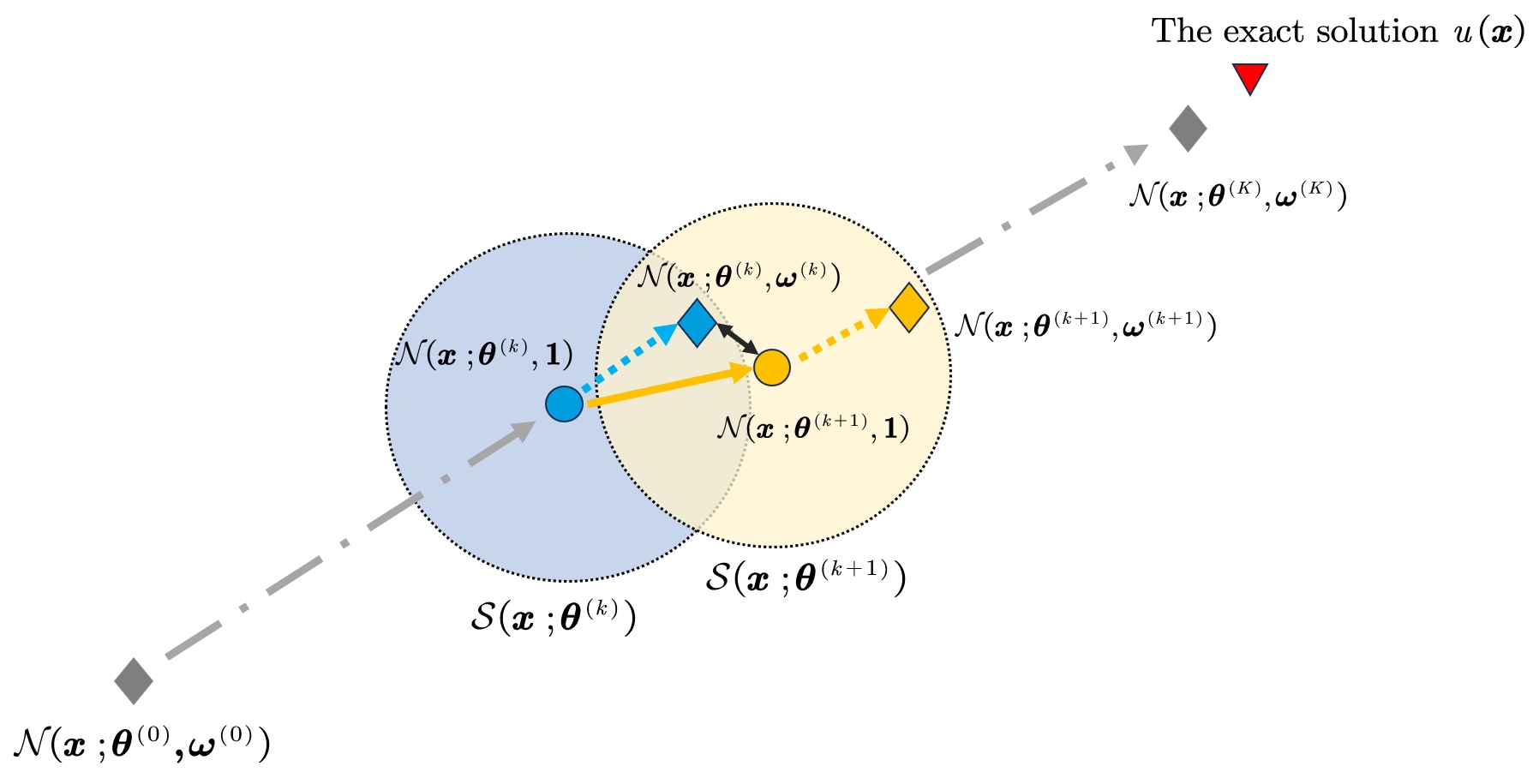}
	\end{center}
	\caption{Diagram of the AO-SNN}
	\label{AOSNN_diagram}
\end{figure}

In this subsection, we explain how the AO-SNN method alternately optimizes the parameters $\bm{\theta}$ and $\bm{\omega}$ and obtains the approximate solution.
In Figure \ref{AOSNN_diagram}, the red triangle in the figure represents the exact solution, the quadrilaterals represent the approximate solutions, the small opaque discs represent the intermediate solutions, and the large translucent discs represent the subspaces. 
The dotted arrows represent the omitted alternating optimization steps, the solid arrows represent the process of training the subspace, and the dashed arrows represent the process of solving the least squares problem.
Assume that the initial approximate solution of the AO-SNN algorithm is $\mathcal{N}(x;\bm{\theta}^{(0)}, \bm{\omega}^{(0)})$ and that after $k~(k=0, 1, \cdots, K-1)$ alternating optimizations we obtain the subspace $\mathcal{S}(x;\bm{\theta}^{(k)})$, the intermediate solution $\mathcal{N}(x;\bm{\theta}^{(k)}, \bm{1})$ and the approximate solution $\mathcal{N}(x;\bm{\theta}^{(k+1)}, \bm{\omega}^{(k+1)})$. 

Thus, in the $(k+1)$th alternating optimization step, we first use the mixed loss function to train the parameters $\bm{\theta}$.
In this way, the intermediate solution $\mathcal{N}(x;\bm{\theta}^{(k)}, \bm{1})$ will be updated to $\mathcal{N}(x;\bm{\theta}^{(k+1)}, \bm{1})$ and the subspace $\mathcal{S}(x;\bm{\theta}^{(k)})$ will be updated to $\mathcal{S}(x;\bm{\theta}^{(k+1)})$.
Affected by the mixed loss function, the intermediate solution $\mathcal{N}(x;\bm{\theta}^{(k+1)}, \bm{1})$ obtained from the new training will not differ much from the previous approximate solution $\mathcal{N}(x;\bm{\theta}^{(k)}, \bm{\omega}^{(k)})$. 
Moreover, since the approximate solution $\mathcal{N}(x;\bm{\theta}^{(k)}, \bm{\omega}^{(k)})$ obtained by the least squares method generally has higher accuracy than the intermediate solution $\mathcal{N}(x;\bm{\theta}^{(k)}, \bm{1})$, the new intermediate solution will be better than the previous one. 
Thus, the subspace $\mathcal{S}(x;\bm{\theta}^{(k)})$ is also effectively optimized into $\mathcal{S}(x;\bm{\theta}^{(k+1)})$. 

Subsequently, we use the least squares method to find a new approximate solution $\mathcal{N}(x;\bm{\theta}^{(k+1)}, \bm{\omega}^{(k+1)})$ in the newer and better subspace $\mathcal{S}(x;\bm{\theta}^{(k+1)})$. 
Repeating these steps, we can obtain the final approximate solution $\mathcal{N}(x;\bm{\theta}^{(K)}, \bm{\omega}^{(K)})$.

\subsection{The implementation of AO-SNN for scattering problems}
In this subsection, we give the implementation details for solving the boundary value problem (\ref{ASwithTTBC}) using the AO-SNN method.
To begin with, in order to approximate the complex-valued solution to (\ref{ASwithTTBC}), it is necessary to design two networks $\mathcal{N}_R(\bm{x};\bm{\theta}_{\rm Re}, \bm{\omega}_{\rm Re})$ and $\mathcal{N}_I(\bm{x};\bm{\theta}_{\rm Im}, \bm{\omega}_{\rm Im})$ for the real part and the imaginary part of the approximate solution, respectively. 
Let $\Phi_{\rm Re}(\bm{x};\bm{\theta}_{\rm Re})$ and $\Phi_{\rm Im}(\bm{x};\bm{\theta}_{\rm Im})$ be the vectors consisting of basis functions in real and imaginary subspace layers, respectively.
Thus, the approximate solution of the AO-SNN method can be expressed as
\begin{equation}\label{cSNNfunction}
	\begin{aligned}
		\mathcal{N}(\bm{x}; \bm{\theta}_{\rm Re}, \bm{\theta}_{\rm Im}, \bm{\omega}_{\rm Re}, \bm{\omega}_{\rm Im}) &= \bm{\omega}_{\rm Re}\cdot\Phi_{\rm Re}(\bm{x};\bm{\theta}_{\rm Re}) + \mathrm{i} \bm{\omega}_{\rm Im}\cdot\Phi_{\rm Im}(\bm{x};\bm{\theta}_{\rm Im}),\\
		&\triangleq\mathcal{N}_{\rm R}(\bm{x};\bm{\theta}_{\rm Re},\bm{\omega}_{\rm Re}) + \mathrm{i}\mathcal{N}_{\rm I}(\bm{x};\bm{\theta}_{\rm Im},\bm{\omega}_{\rm Im}).
	\end{aligned}
\end{equation}

To solve the boundary value problem (\ref{ASwithTTBC}) using the AO-SNN method, we need to optimize the parameters $(\bm{\theta}_{\rm Re}, \bm{\theta}_{\rm Im})$ and $(\bm{\omega}_{\rm Re}, \bm{\omega}_{\rm Im})$ alternately as described in Section 3.2. We also need to set $\bm{\omega}_{\rm Re}=\bm{\omega}_{\rm Im}=\bm{1}$, initialize $\bm{\theta}_{\rm Re}^{(k)}=\bm{\theta}_{\rm Re}^{(k-1)}$ and $ \bm{\theta}_{\rm Im}^{(k)}=\bm{\theta}_{\rm Im}^{(k-1)}$, and use a mixed loss function for each alternating optimization.
Denote by $X_{\Omega}$, $X_{\partial D}$, and $X_{TBC}$ the sets of collocation points in $\Omega$, on $\partial D$, and on $\partial B_R$, respectively.
The mixed loss function for the $k$th iteration can be defined as
\begin{align*}
	&\mathcal{L}^{\rm Mixed}\left(\mathcal{N}(\bm{x}; \bm{\theta}_{\rm Re}^{(k)}, \bm{\theta}_{\rm Im}^{(k)}, \bm{1}, \bm{1})\right) \\
	&=  \eta_k\mathcal{L}^{\rm Metric}\left(\mathcal{N}(\bm{x}; \bm{\theta}_{\rm Re}^{(k)}, \bm{\theta}_{\rm Im}^{(k)}, \bm{1}, \bm{1}), \mathcal{N}(\bm{x}; \bm{\theta}_{\rm Re}^{(k-1)}, \bm{\theta}_{\rm Im}^{(k-1)}, \bm{\omega}_{\rm Re}^{(k-1)}, \bm{\omega}_{\rm Im}^{(k-1)})\right) \\
	&\quad +  \sigma_k\mathcal{L}^*\left(\mathcal{N}(\bm{x}; \bm{\theta}_{\rm Re}^{(k)}, \bm{\theta}_{\rm Im}^{(k)}, \bm{1}, \bm{1})\right),
\end{align*}
where
\begin{equation}
	\begin{aligned}
		&\mathcal{L}^{\rm Metric}\left(\mathcal{N}(\bm{x}; \bm{\theta}_{\rm Re}^{(k)}, \bm{\theta}_{\rm Im}^{(k)}, \bm{1}, \bm{1}), \mathcal{N}(\bm{x}; \bm{\theta}_{\rm Re}^{(k-1)}, \bm{\theta}_{\rm Im}^{(k-1)}, \bm{\omega}_{\rm Re}^{(k-1)}, \bm{\omega}_{\rm Im}^{(k-1)})\right) \\
		&=\frac{1}{N}\sum^{N}_{i=1}\sum_{|\bm{\alpha}|=0}^{\gamma}\big(\partial^{\bm{\alpha}}\mathcal{N}_{\rm R}(\bm{x}_i; \bm{\theta}_{\rm Re}^{(k)}, \bm{1})-\partial^{\bm{\alpha}}\mathcal{N}_{\rm R}(\bm{x}_i; \bm{\theta}_{\rm Re}^{(k-1)}, \bm{\omega}_{\rm Re}^{(k-1)}\big)^2\\
		& \quad+\frac{1}{N}\sum^{N}_{i=1}\sum_{|\bm{\alpha}|=0}^{\gamma}\big(\partial^{\bm{\alpha}}\mathcal{N}_{\rm I}(\bm{x}_i; \bm{\theta}_{\rm Im}^{(k)}, \bm{1})-\partial^{\bm{\alpha}}\mathcal{N}_{\rm I}(\bm{x}_i; \bm{\theta}_{\rm Im}^{(k-1)}, \bm{\omega}_{\rm Im}^{(k-1)})\big)^2,
	\end{aligned}
\end{equation}
and 
\begin{equation}
	\begin{aligned}
		&\mathcal{L}^*(\mathcal{N}(\bm{x})) \\
		&= \frac{1}{N_{\rm \Omega}}\sum_{\bm{x}\in X_{\rm \Omega}}\Big(\big(\Delta\mathcal{N}_{\rm R}(\bm{x}) + \kappa^2 \mathcal{N}_{\rm R}(\bm{x})\big)^2 + (\Delta\mathcal{N}_{\rm I}(\bm{x}) + \kappa^2 \mathcal{N}_{\rm I}(\bm{x}))^2\Big)\\
		&\quad+  \frac{1}{N_{\rm \partial D}}\sum_{\bm{x}\in X_{\rm \partial D}}\Big(\big(\mathcal{N}_{\rm R}(\bm{x}) - {\rm Re}\{g(\bm{x})\}\big)^2 + \big(\mathcal{N}_{\rm I}(\bm{x}) - {\rm Im}\{g(\bm{x})\}\big)^2\Big) \\
		&\quad+ \frac{1}{N_{\rm TBC}}\sum_{\bm{x}\in X_{\rm TBC}}\left(\Big(\frac{\partial\mathcal{N_{\rm R}}(\bm{x})}{\partial r} - {\rm Re}\{\mathcal{F}_N\mathcal{N}(\bm{x})\}\Big)^2 + \Big(\frac{\partial\mathcal{N}_{\rm I}(\bm{x})}{\partial r} - {\rm Im}\{\mathcal{F}_N\mathcal{N}(\bm{x})\}\Big)^2\right).
	\end{aligned}
\end{equation}
Here, $N_{\rm \Omega}, N_{\rm \partial D}$ and $N_{\rm TBC}$ represent the numbers of 
collocation points in $\Omega$, on $\partial D$ and on $\partial B_R$, respectively.

Once the neural network is adequately trained, we obtain the parameters $\bm{\theta}_{\rm Re}^{(k)}$ and $\bm{\theta}_{\rm Im}^{(k)}$. 
Since all operators in boundary value problem (\ref{ASwithTTBC}) are linear operators, we can compute the parameters $\bm{\omega}_{\rm Re}^{(k)}$ and $\bm{\omega}_{\rm Im}^{(k)}$ by solving the least square problem
\begin{equation}\label{AS_lstsq}
	\argmin_{(\bm{\omega}_{\rm Re}^{(k)}, \bm{\omega}_{\rm Im}^{(k)})}	\Big(||\Delta \mathcal{N}(\bm{x}) + \kappa^2\mathcal{N}(\bm{x})||_{l^2({X_{\Omega}})}^2 + ||\mathcal{N}(\bm{x}) - g(\bm{x})||_{l^2(X_{\partial D})}^2 + ||\mathcal{F}_N[\mathcal{N}(\bm{x})] - \frac{\partial\mathcal{N}(\bm{x})}{\partial r}||_{l^2(X_{TBC})}^2\Big).
\end{equation}
Thus, we obtain an approximate solution $\mathcal{N}(\bm{x}; \bm{\theta}_{\rm Re}^{(k)}, \bm{\theta}_{\rm Im}^{(k)}, \bm{\omega}_{\rm Re}^{(k)}, \bm{\omega}_{\rm Im}^{(k)})$. 

\section{Numerical Experiments}\label{NE}

In this section, we present several numerical examples to illustrate the performance of AO-SNN method  and to provide comparisons with other commonly used methods in solving acoustic scattering problems (\ref{ASwithTTBC}).

We denote the relative $l^2$ error by
\begin{equation}
	\mathcal{E}_{l^2}\big(\mathcal{N}(\bm{x})\big)=\frac{\sqrt{\sum_{\bm{x}\in \mathcal{X}}|\mathcal{N}(\bm{x}) - u(\bm{x})|^2}}{\sqrt{\sum_{\bm{x}\in \mathcal{X}}|u(\bm{x})|^2}}, 
\end{equation}
where $\mathcal{X}$ represents the set of all collocation points, and $u(\bm{x})$ is the exact solution of the PDEs.

We employ the deep learning framework PyTorch and the least squares algorithm from the SciPy library in or algorithm. 
All data are in double precision format.

\begin{example}
	In this example, the scatter $D$ is defined as the circle of radius $r=0.5$, and the TBC is set on the circle $B_R$ with $R=1$. 
	The exact solution of this example is represented by $H^{(1)}_{0}(\kappa r)$. Unless stated, we take $\kappa = 5$.
\end{example}
\begin{figure}[htp!]
	\begin{minipage}{0.45\linewidth}
		\begin{center}
			\includegraphics[scale=0.5]{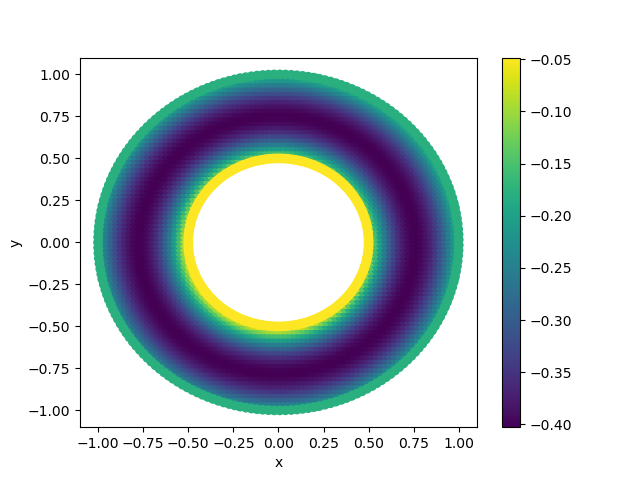}
		\end{center}
	\end{minipage}
	\begin{minipage}{0.45\linewidth}
		\begin{center}
			\includegraphics[scale=0.5]{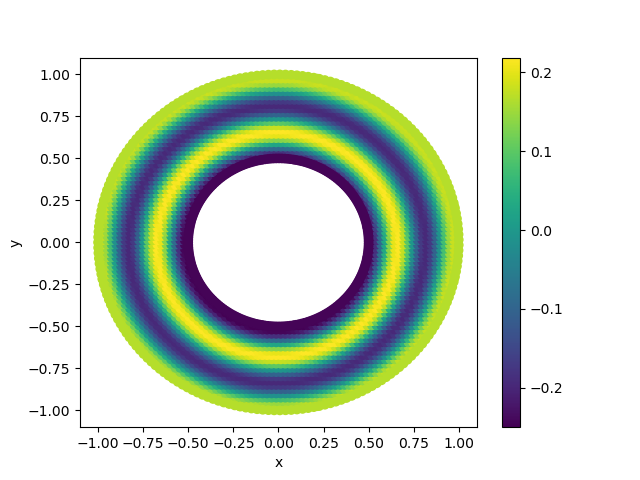}
		\end{center}
	\end{minipage}
	\caption{The real part of exact solution (left: $\kappa=5$, right: $\kappa=20$)}
	\label{example1,fig}
\end{figure}
We first compare the effectiveness of the PINN method \cite{PINN_2019}, the RBDNN method \cite{RBDNN_2022, RBDNN_2023}, the SNN method \cite{SNN_arXiv}  and the AO-SNN method in solving Example 1 with sound soft boundary condition. 
For the SNN method and the AO-SNN method, we set the depth of the hidden layer to 3 with a width of 40 neurons, and set the width of subspace layer to 600 neurons. This makes the total number of parameters 27880. The maximum number of iterations of the AO-SNN method is $K=2$.
In order to ensure a fair comparison in terms of the number of parameters, we set the hidden layer with a depth of 4 and a width of 100 neurons for the PINN method and the RBDNN method. This makes the total number of parameters 30700. 
In all neural networks, we use $\mathrm{tanh}(x)$ as the activation function, and choose the Adam method with the learning rate $0.001$ as the optimizer.
\begin{figure}[htp!]
	\subcaptionbox{PINN}{
		\label{subfig1}
		\includegraphics[scale=0.4]{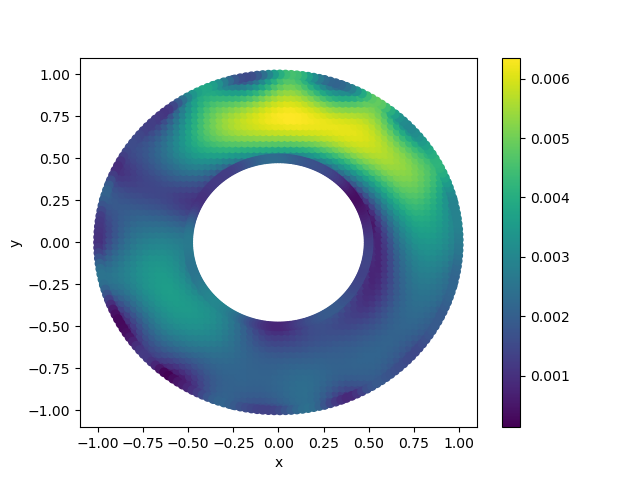
	}}
	\hspace{.4in}
	\subcaptionbox{RBDNN}{
		\label{subfig2}
		\includegraphics[scale=0.4]{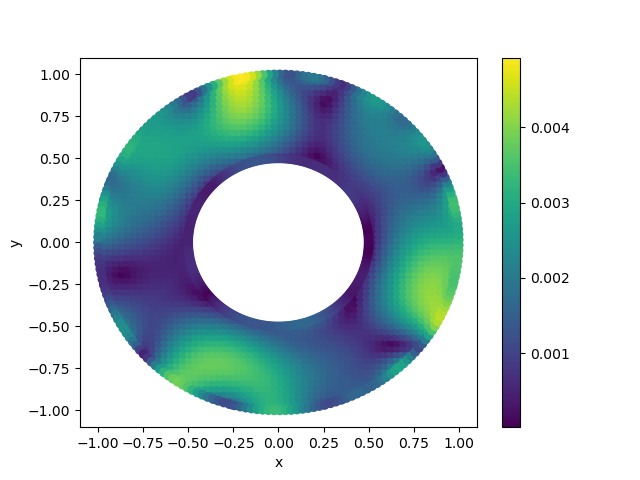
	}}
	\newline
	\subcaptionbox{SNN}{
		\label{subfig3}
		\includegraphics[scale=0.4]{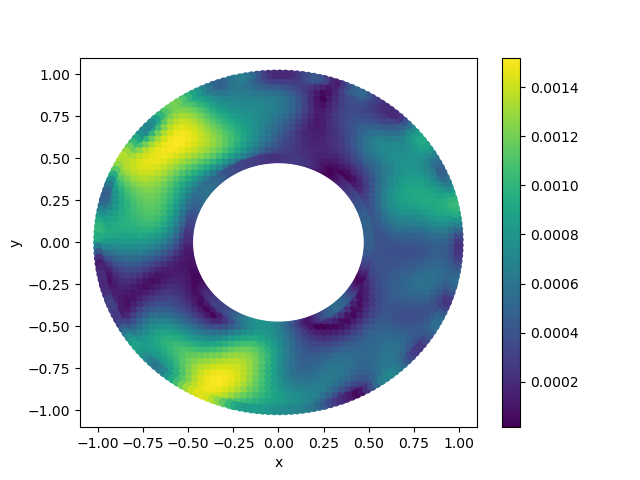
	}}
	\hspace{.4in}
	\subcaptionbox{AO-SNN}{
		\label{subfig4}
		\includegraphics[scale=0.4]{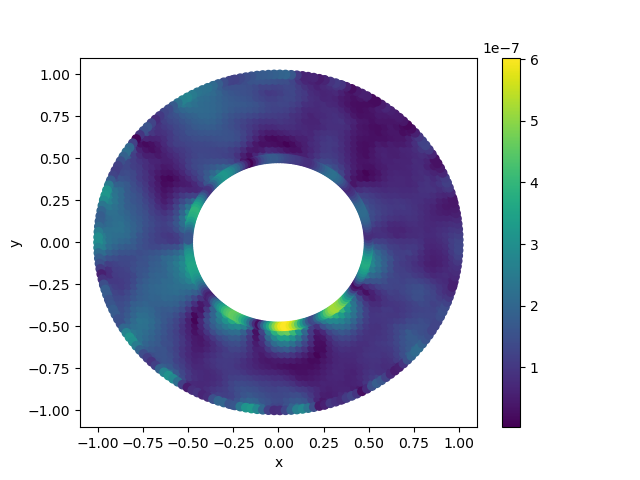
	}}
	\caption{The pointwise errors of PINN, RBDNN, SNN and AO-SNN}
	\label{SNN,PINN,H0,fig}
\end{figure}

\begin{table}[!htbp]
	\caption{Errors and epochs of PINN, RBDNN, SNN, and AO-SNN}
	\label{SNN,PINN,H0,table}
	\centering
	\begin{tabular}{c||c|c}
		\toprule[2pt]
		Name & Relative $l^2$ error & Epochs  \\\midrule
		PINN & $7.381\times 10^{-3}$ & 40012 \\   
		RBDNN & $5.208\times10^{-3}$ & 40001 \\   
		SNN & $1.731\times10^{-3}$ & 14000 \\   
		AO-SNN & $3.693\times10^{-7}$ & 9006  \\   
		\bottomrule[2pt]
	\end{tabular}
\end{table}
Table \ref{SNN,PINN,H0,table} lists the relative $l^2$ errors of the PINN, RBDNN, SNN and AO-SNN methods, and the numbers of epochs required for each method. 
The methods selected here require almost the same time in every epoch because they all need to calculate the second order derivatives in every epoch and have a similar number of parameters. So the number of epochs can also represent the computation time. Though the standard SNN method has the potential to achieve high accuracy in some model problems \cite{SNN_arXiv}, it can not get high accuracy for scattering problems. 
In Figure \ref{SNN,PINN,H0,fig}, the pointwise errors of these methods is illustrated.
It is easy to see that the AO-SNN exhibits superior accuracy compared to the PINN, the RBDNN and the SNN.

\begin{table}[!htbp]
	\caption{Errors and epochs of AO-SNN with different $K$}
	\label{AOSNN,iteration,H0,table}
	\centering
	\begin{tabular}{c||c|c}
		\toprule[2pt]
		$K$ & Relative $l^2$ error & Epochs \\\midrule
		0 & $8.184\times10^{-3}$ & 1000 \\
		1 & $1.265\times10^{-5}$ & 3250 \\   
		2 & $3.582\times 10^{-6}$ & 4906 \\   
		3 & $4.804\times10^{-7}$ & 11722 \\   
		4 & $1.591\times10^{-8}$ & 61722  \\   
		5 & $5.170\times10^{-8}$ & 111722  \\   
		\bottomrule[2pt]
	\end{tabular}
\end{table}

Table \ref{AOSNN,iteration,H0,table} shows how the relative $l^2$ error decays as the number of alternating optimization increases. 
The training of the subspace is terminated when $\mathcal{L}_{Mixed}$ reaches 0.1 times of itself at the beginning, or when the number of epoch reaches $50000$.
It is easy to see when $K=4$, the AO-SNN reaches the accuracy of $10^{-8}$.
When $K=4$ and $5$, the number of epochs in the last alternating optimization step reaches $50000$.
The reason for the lack of further improvement in accuracy at $K=5$ may be due to some error in the truncation of the DtN operator in the boundary value problem (\ref{ASwithTTBC}).

In Table \ref{AOSNN,MN,H0,table}, we shows the errors of the AO-SNN with different the number of collocation points $N$ and the dimension of the subspace $M$.
When $M$ is fixed, the increase of $N$ has little effect on the improvement of accuracy, while when $N$ is fixed, the increase of $M$ has a significant effect on the improvement of accuracy. 
This indicates that the dimension of the subspace has a greater impact on accuracy than the number of collocation points.
However, the larger dimension of the subspace means that the least squares stage is more computationally expensive.

\begin{table}[!htbp]
	\caption{Errors of AO-SNN with different $N$ and $M$}
	\label{AOSNN,MN,H0,table}
	\centering
	\begin{tabular}{c||c|c|c|c}
		\toprule[2pt]
		\diagbox{$N$}{$M$} & 50 & 100 & 300 & 600  \\ \midrule
		606 & $3.904\times10^{-2}$ & $2.449\times10^{-3}$ & $1.824\times10^{-3}$ & $9.342\times10^{-4}$\\   
		1282 & $5.757\times10^{-2}$ & $1.684\times10^{-3}$ & $1.123\times10^{-5}$ & $7.884\times10^{-6}$\\   
		2434 & $4.127\times10^{-2}$ & $3.871\times10^{-3}$ & $2.747\times10^{-6}$ & $3.808\times10^{-7}$\\   
		4046 & $1.264\times10^{-1}$ & $5.772\times10^{-3}$ & $2.954\times10^{-6}$ & $1.457\times10^{-7}$\\
		\bottomrule[2pt]
	\end{tabular}
\end{table}

Table \ref{AOSNN,Layers,H0,table} demonstrates how the depth of the hidden layer affects the accuracy of the AO-SNN.
In this test, we set $K=1$.
Obviously, the numerical accuracy of the AO-SNN is very low if there is no hidden layer, while the accuracy can reach $10^{-6}$ with only one hidden layer.
The greater the depth of the hidden layer, the smaller the number of epochs required in the AO-SNN.

\begin{table}[!htbp]
	\caption{Errors and epochs of AO-SNN with different number of hidden layers}
	\label{AOSNN,Layers,H0,table}
	\centering
	\begin{tabular}{c||c|c}
		\toprule[2pt]
		Number of hidden layers & Relative $l^2$ error & Epochs  \\ \midrule
		0 & $3.513\times10^{-3}$ & 25000 \\
		1 & $8.242\times10^{-6}$ & 17173 \\ 
		2 & $8.663\times10^{-7}$ & 15257 \\   
		3 & $4.130\times10^{-7}$ & 8590  \\
		4 & $3.005\times10^{-7}$ & 6872 \\
		5 & $2.432\times10^{-7}$ & 6315 \\	   
		\bottomrule[2pt]
	\end{tabular}
\end{table}

Table \ref{SNN,orders,H0,table} compares accuracy of AO-SNN method with different $\gamma$ in mixed loss function (\ref{Mixed_loss}).
We define 
\begin{equation}\label{moderror}
	\mathcal{E}_{H^{\beta}} = \sqrt{\sum_{x\in\mathcal{X}}\sum_{|\alpha|=\beta}\big(\partial^{\alpha}\mathcal{N}(x)-\partial^{\alpha}u(x)\big)^2}, \quad \beta = 0,1,2
\end{equation}
to test the accuracy of the AO-SNN solution.

When $\gamma \le2$, it can be observed that the larger the value of $\gamma$, the higher the accuracy. However, when $\gamma=3$, the accuracy does not improve further. This phenomenon indicates that a good choice for $\gamma$ is the highest order of partial derivatives in the target equation, i.e., $2$ for problem \ref{ASwithTTBC}.

\begin{table}[!htbp]
	\caption{Errors of AO-SNN in different norms}
	\label{SNN,orders,H0,table}
	\centering
	\begin{tabular}{c||c|c|c|c}
		\toprule[2pt]
		$\gamma$ & Relative $l^2$ error & $\mathcal{E}_{H^{0}}$ & $\mathcal{E}_{H^{1}}$ & $\mathcal{E}_{H^{2}}$   \\ \midrule
		0 & $2.028\times10^{-5}$ & $1.289\times10^{-5}$ & $8.191\times10^{-5}$ & $1.218\times10^{-3}$  \\     
		1 & $7.532\times10^{-6}$ & $4.788\times10^{-6}$ & $4.638\times10^{-5}$ & $1.246\times10^{-3}$ \\   
		2 & $7.640\times10^{-8}$ & $4.857\times10^{-8}$ & $5.468\times10^{-7}$ & $1.556\times10^{-5}$  \\   
		3 & $8.667\times10^{-8}$ & $5.509\times10^{-8}$ & $7.782\times10^{-7}$ & $2.870\times10^{-5}$ \\
		\bottomrule[2pt]
	\end{tabular}
\end{table}

In Table \ref{SNN,BCs,H0,table}, we show the numerical solution of the AO-SNN method in solving problem (\ref{ASwithTTBC}) with various boundary conditions.
We take $\lambda=1$ in impedance boundary condition, and set $K=1$. The accuracy of our algorithm is not sensitive to boundary conditions.

\begin{table}[!htbp]
	\caption{Errors and epochs of AO-SNN for scattering problem with different boundary conditions}
	\label{SNN,BCs,H0,table}
	\centering
	\begin{tabular}{c||c|c}
		\toprule[2pt]
		Boundary condition & Relative $l^2$ error & Epochs  \\ \midrule
		Sound-soft & $7.947\times10^{-7}$ & 11591  \\   
		Sound-hard & $5.509\times10^{-7}$ & 13950 \\   
		Impedance & $4.425\times10^{-7}$ & 15138 \\     
		\bottomrule[2pt]
	\end{tabular}
\end{table}

We compare the AO-SNN method with several other methods in solving scattering problems with large wavenumbers.
Although the AO-SNN method cannot achieve $10^{-7}$ accuracy at wavenumber $\kappa=20$, it can still achieve $10^{-4}$ accuracy and outperforms several other classes of methods; see Table \ref{SNN,k20,H0,table}.
In Figure \ref{SNN,PINN,k20,figure}, we show the real parts of the exact, PINN, SNN and AO-SNN solutions.

\begin{table}[!htbp]
	\caption{Errors and epochs of  PINN, RBDNN, SNN and AO-SNN for $\kappa=20$}
	\label{SNN,k20,H0,table}
	\centering
	\begin{tabular}{c||c|c}
		\toprule[2pt]
		Name & Relative $L^2$ error & Epochs  \\\midrule
		PINN & $0.999$ & 50000 \\  
		RBDNN & $1.450$ & 50000  \\ 
		SNN & $0.566$ & 46000 \\
		AO-SNN & ${9.126\times10^{-4}}$ & 28833  \\
		\bottomrule[2pt]
	\end{tabular}
\end{table}
\begin{figure}[htp!]
	\subcaptionbox{Exact solution}{
		\label{Exact solution}
		\includegraphics[scale=0.4]{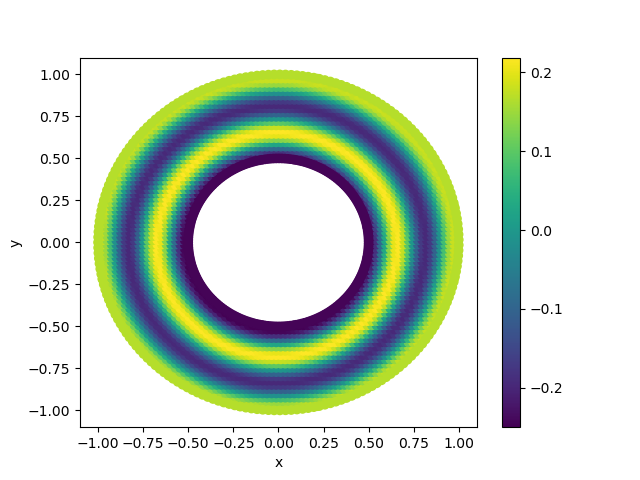
	}}
	\hspace{.4in}
	\subcaptionbox{PINN solution}{
		\label{subfig2}
		\includegraphics[scale=0.4]{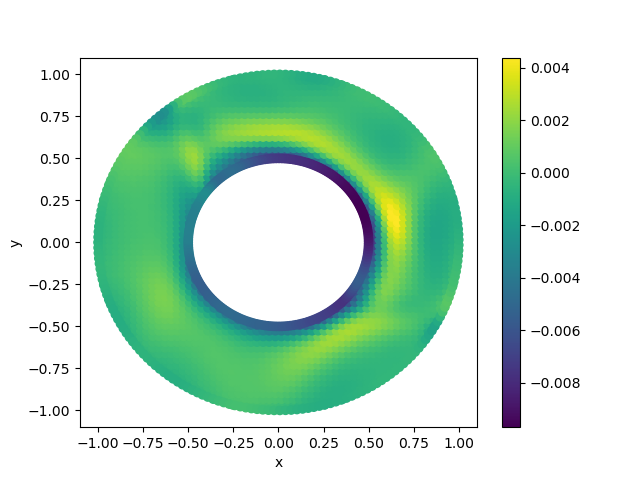
	}}
	\newline
	\subcaptionbox{SNN solution}{
		\label{subfig3}
		\includegraphics[scale=0.4]{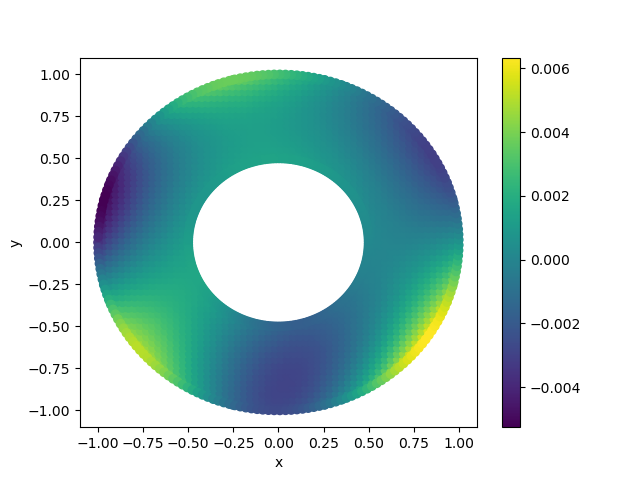
	}}
	\hspace{.4in}
	\subcaptionbox{AO-SNN solution}{
		\label{subfig4}
		\includegraphics[scale=0.4]{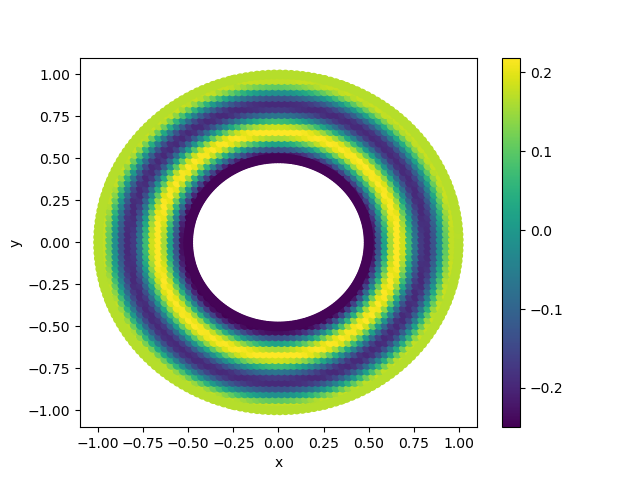
	}}
	\caption{The real part of numerical solutions of PINN, SNN and AO-SNN for scattering problem with $\kappa=20$}
	\label{SNN,PINN,k20,figure}
\end{figure}

\begin{example}
	In this example, the scatterer $D$ is defined as a circle of radius $r=0.5$, and the TBC is set on a circle $B_R$ with $R=1$. 
	The incident wave is a plane wave propagating in the direction $(1,0)$, with the wavenumber $\kappa=5$. 
\end{example}

We set the hidden layer depth of SNN and AO-SNN to 3, the width of each hidden layer to 40, and the width of the subspace layer to 1000.
This means that the total number of parameters for both methods is 44280. 
The PINN network is configured as 4 hidden layers with 125 neurons in each layer, resulting in a total of 47,250 parameters.
We choose $K=2$ for the AO-SNN, $\tanh(x)$ as activation function for all networks, and the Adam method as optimizer with the learning rate 0.001.

In Table \ref{SNN,PINN,inwave,table}, we show numerical results of the PINN, the SNN and the AO-SNN with different hyperparameters.
Among these methods only the AO-SNN is able to obtain numerical results with high accuracy.
In the mixed loss function $\mathcal{L}^{\rm Mixed}$, the numerical results corresponding to $\sigma_k=0$ are better than those corresponding to $\sigma_k=1$. This indeed indicates that the metric loss function $\mathcal{L}^{\rm Metric}$ occupies a more important place in the mixed loss function.
Figure \ref{SNN,PINN,inwave,figure} illustrates the pointwise errors of the PINN, the SNN and the AO-SNN with different hyperparameters.

\begin{table}[!htbp]
	\caption{The relative $l^2$ errors and epochs of PINN, SNN and AO-SNN with different hyperparameters}
	\label{SNN,PINN,inwave,table}
	\centering
	\begin{tabular}{c||c|c}
		\toprule[2pt]
		Name & Relative $l^2$ error & Epochs  \\\midrule
		PINN & $6.624\times 10^{-3}$ & 50000 \\ 
		SNN & $4.750\times10^{-2}$ & 46000 \\ 
		AO-SNN ($\eta_k=1, \sigma_k=1$) & $1.829\times10^{-5}$ & 41000  \\ 
		AO-SNN ($\eta_k=1, \sigma_k=0$) & $1.204\times10^{-6}$ & 28833  \\
		\bottomrule[2pt]
	\end{tabular}
\end{table}

\begin{figure}[htp!]
	\subcaptionbox{PINN}{
		\label{subfig1}
		\includegraphics[scale=0.4]{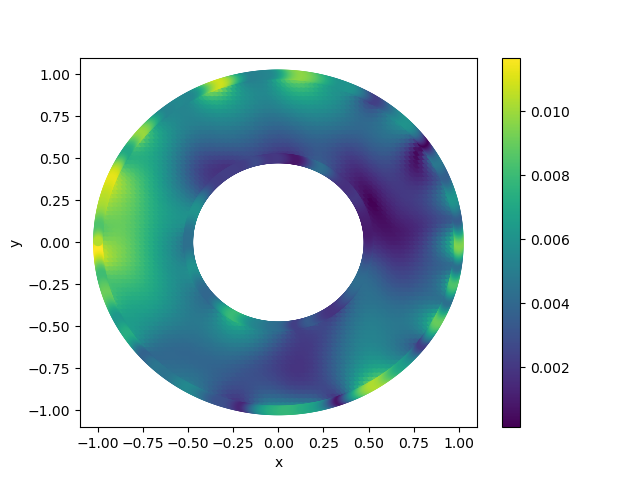
	}}
	\hspace{.4in}
	\subcaptionbox{SNN}{
		\label{subfig2}
		\includegraphics[scale=0.4]{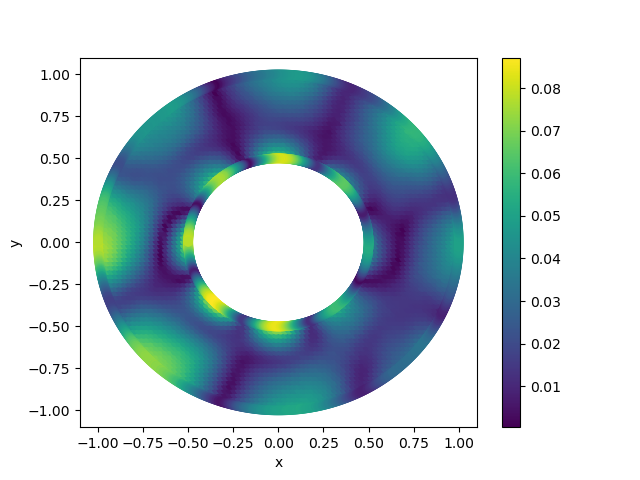
	}}
	\newline
	\subcaptionbox{AO-SNN ($\eta_k=1, \sigma_k=1$)}{
		\label{subfig3}
		\includegraphics[scale=0.4]{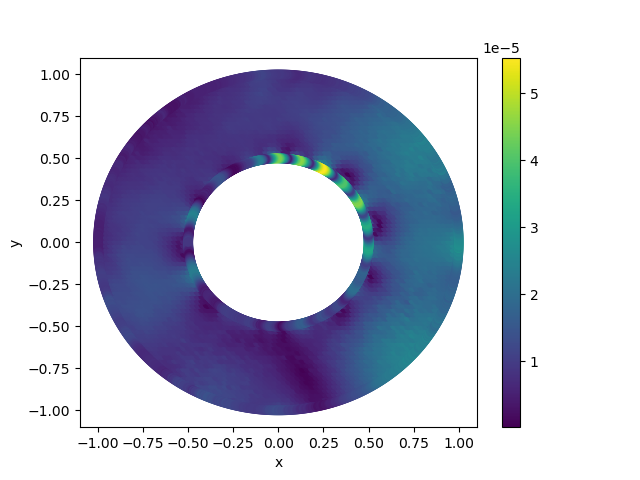
	}}
	\hspace{.4in}
	\subcaptionbox{AO-SNN ($\eta_k=1, \sigma_k=0$)}{
		\label{subfig4}
		\includegraphics[scale=0.4]{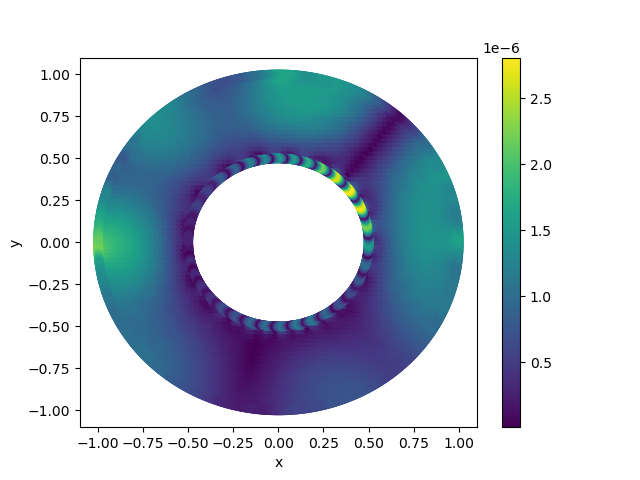
	}}
	\caption{The pointwise errors of PINN, SNN, and AO-SNN with different hyperparameters}
	\label{SNN,PINN,inwave,figure}
\end{figure}

\section{Discussions}
In this section, we try to present some discussions on the accuracy of the AO-SNN method.
Let $\Pi$ represent the projection operator that maps any sufficient smooth function $u$ into the neural network subspace $\mathcal{S}$. Thus, we have
\begin{equation}
	\begin{aligned}
		\mathcal{E}^2_{l^2}\big(\mathcal{N}(\bm{x})\big)&=\frac{\sum_{\bm{x}\in \mathcal{X}}| \mathcal{N}(\bm{x}) - u(\bm{x}) |^2}{\sum_{x\in \mathcal{X}}| u(\bm{x}) |^2},\\
		&\le\frac{\sum_{\bm{x}\in \mathcal{X}}| \mathcal{N}(\bm{x}) - \Pi u(\bm{x}) |^2 + \sum_{\bm{x}\in \mathcal{X}}| \Pi u(\bm{x}) - u(\bm{x}) |^2}{\sum_{\bm{x}\in \mathcal{X}}| u(\bm{x}) |^2},\\
		&=\frac{\sum_{\bm{x}\in \mathcal{X}}|\mathcal{N}(\bm{x}) - \Pi u(\bm{x})|^2}{\sum_{\bm{x}\in \mathcal{X}}| u(\bm{x}) |^2} + \frac{\sum_{\bm{x}\in \mathcal{X}}| \Pi u(\bm{x}) - u(\bm{x}) |^2}{\sum_{\bm{x}\in \mathcal{X}}| u(\bm{x}) |^2},\\
		&\triangleq \mathcal{E}_{ls}(\mathcal{N}(\bm{x})) + \mathcal{E}_{space}(\mathcal{N}(\bm{x})),
	\end{aligned}
\end{equation}
where $\mathcal{E}_{ls}$ represents the error of the least squares solution, and $\mathcal{E}_{space}$ represents the projection error.

In fact, a starting point for the design of the AO-SNN method is to find the proper subspace making the projection error $\mathcal{E}_{space}$ very small.
Once the projection error is sufficiently small, the error of the AO-SNN solution is determined by the least squares error $\mathcal{E}_{ls}$.

The process of solving optimization problems using least squares (17) is equivalent to finding a least squares solution to a system of hyper-determined linear equations
\begin{equation}\label{linear_eq}
	A\bm{\omega} = \bm{b},
\end{equation}
where 
\begin{equation}\label{MatrixA}
	A=\begin{bmatrix}
		A^{pde}\\
		A^{bc}
	\end{bmatrix},~~A^{pde}_{i,j} = \phi_j(x_i),~~A^{bc}_{i,j} = \phi_j(\bar{x}_i),
\end{equation}
and 
\begin{equation}
	\bm{b}=\begin{bmatrix}
		\bm{b}^{pde}\\
		\bm{b}^{bc}
	\end{bmatrix},~~b^{pde}_{i} = f(x_i),~~b^{bc}_{i} = g(\bar{x}_i).
\end{equation}
As illustrated in equation (\ref{MatrixA}), the matrix $A$ has $N+\bar{N}$ rows and $M$ columns. Since the number is generally much larger than the number of basis functions, equation (\ref{linear_eq}) is hyper-determined and has at least one least squares solution.

Denote the $l^2$ vector norm by $\left\|\cdot\right\|_{l^2}$, and the condition number of matrix $A$ by $C_A=\left\|A\right\|_{l^2}\cdot\left\|A\right\|_{l^2}^{-1}$. Suppose that the right hand term $\bm{b}$ is perturbed by $\delta \bm{b}$, such that the numerical solution changes $\delta \bm{\omega}$. Therefore, the norm of the error $\delta \bm{\omega}$ has the following estimate
\begin{equation}\label{LS_estimate}
	{C_A}^{-1}\left\| \bm{\omega} \right\|_{l^2} \frac{\left\| \delta \bm{b} \right\|_{l^2}}{\left\| \bm{b} \right\|_{l^2}}\,\,\le \,\,\left\| \delta \bm{\omega} \right\|_{l^2} \,\,\le \,\,C_A\left\| \bm{\omega} \right\|_{l^2} \frac{\left\| \delta \bm{b} \right\|_{l^2}}{\left\| \bm{b} \right\|_{l^2}}.
\end{equation}
Since the numerical solution is obtained by a linear combination of basis functions; see (\ref{forward_SNN}), the estimate (\ref{LS_estimate}) ultimately results in a relative error of the numerical solution

\begin{equation}\label{Els}
	\begin{aligned}
		\mathcal{E}^2_{ls}\big(\mathcal{N}(\bm{x})\big)
		&=\frac{\sum_{\bm{x}\in \mathcal{X}}|\sum_{j=1}^{M}\delta\omega_j\phi_j(\bm{x}) |^2}{\sum_{\bm{x}\in \mathcal{X}}| u(\bm{x}) |^2}\\
		&= \frac{\sum_{\bm{x}\in \mathcal{X}}(\delta\bm{\omega}\cdot{\Phi}(\bm{x}))^2}{\sum_{\bm{x}\in \mathcal{X}}| u(\bm{x}) |^2}\\
		&\le \frac{\left\|\delta\bm{\omega}\right\|_{l^2}^2\sum_{\bm{x}\in \mathcal{X}}(\left\|\Phi(\bm{x})\right\|_{l^2}^2)}{\sum_{\bm{x}\in \mathcal{X}}| u(\bm{x}) |^2}\\
		&\le \frac{\sum_{\bm{x}\in \mathcal{X}}(\left\|\Phi(\bm{x})\right\|_{l^2}^2)}{\sum_{\bm{x}\in \mathcal{X}}| u(\bm{x}) |^2} C_A^2\left\| \bm{\omega} \right\|_{l^2}^2 \frac{\left\| \delta \bm{b} \right\|_{l^2}^2}{\left\| \bm{b} \right\|_{l^2}^2},
	\end{aligned}
\end{equation}
where $u$ represents the exact solution of the PDEs. It follows from (\ref{Els}) that the error $\mathcal{E}^2_{ls}$ can be controlled by $\left\| \bm{\omega} \right\|^2_{l^2}$.
In Figure \ref{LS}, we present the errors and the corresponding norms $\left\|\bm{\omega}\right\|_{l^2}$ for the AO-SNN solution for different alternating optimization step $K$'s.
It is clear that as $K$ increases, the rate at which the exponent of $\left\|\bm{\omega}\right\|_{l^2}$ becomes smaller is almost the same as the rate at which the exponent of relative $l^2$ error becomes smaller.

\begin{table}[!htbp]
	\caption{The error and $\left\|\bm{\omega}\right\|_{l^2}$ for different  alternating optimization step $K$'s}
	\label{LS}
	\centering
	\begin{tabular}{c||c|c}
		\toprule[2pt]
		$K$ & Relative $l^2$ error& $\left\|\bm{\omega}\right\|_{l^2}$  \\\midrule
		0 & $1.2\times10^{-3}$ & $7.321\times10^{10}$  \\ 
		1 & $4.2680\times10^{-5}$ & $3.304\times10^{9}$  \\ 
		3 & $1.6740\times10^{-6}$ & $1.742\times10^{8}$  \\ 
		4 & $3.9967\times10^{-7}$ &$ 8.359\times10^{7}$ \\
		5 & $8.0010\times10^{-7}$ & $7.913\times10^{7}$  \\
		6 & $8.6645\times10^{-7}$ & $6.170\times10^{7}$   \\
		\bottomrule[2pt]
	\end{tabular}
\end{table}

\section{Conclusion}
In this paper, we reformulate the acoustic scattering problem in unbounded domain into a boundary value problem in bounded domain by the TBC method. We introduce a high-accuracy machine learning-based method, the alternately-optimized subspace method based on neural network (AO-SNN), to solve this boundary value problem.

In the AO-SNN method, our innovation is twofold: firstly, we give a new seemingly interpretable way of updating the subspace; secondly, we give a new loss function to measure the distance between the approximate solution of the previous step and the new subspace.
Our proposed AO-SNN method has the highest accuracy among the various machine learning-based methods currently available for solving scattering problems in unbounded domains.

However, the accuracy of AO-SNN in solving scattering problems in unbounded domains does not reach the accuracy of the SNN method in solving the usual boundary value problems. 
The reason may be that the numerical computation of the scattering problem contains some truncation of the infinite series, such as the truncation of the DtN operator.

In future work, we hope to further improve the accuracy of AO-SNN for solving scattering problems in unbounded regions. In addition, we will explore the feasibility of the AO-SNN method for solving large wavenumber problems.

\section*{Acknowledgements}
This work of Z.S. is partially supported by the Fund of the National Key Laboratory of Computational Physics grant 6142A05230501. D.W. is partially supported by the National Natural Science Foundation of China grant 12422116, Guangdong Basic and Applied Basic Research Foundation grant 2023A1515012199, Shenzhen Science and Technology Innovation Program grants JCYJ20220530143803007 and RCYX20221008092843046, and Hetao Shenzhen-Hong Kong Science and Technology Innovation Cooperation Zone Project grant HZQSWS-KCCYB-2024016. This work of J.L. is partially supported by the National Natural Science Foundation of China grant 12271209.

\section*{Declarations}

\noindent
\textbf{Conﬂict of interest} The authors declare that they have no conﬂict of interest.
\vskip 5pt
\noindent
\textbf{Data availability} The data that support the ﬁndings of this study are available from the corresponding author upon reasonable request.


\begin{thebibliography}{00}

\bibitem{PML_1994}Berenger J P. A perfectly matched layer for the absorption of electromagnetic waves. Journal of computational physics, 1994, 114(2): 185-200.

\bibitem{PhaseDNN_2020}Cai W, Li X, Liu L. A phase shift deep neural network for high frequency approximation and wave problems. SIAM Journal on Scientific Computing, 2020, 42(5): A3285-A3312.

\bibitem{DFVM_2024}Cen J, Zou Q. Deep finite volume method for partial differential equations. Journal of Computational Physics, 2024, 517: 113307.

\bibitem{RFM_2022}Chen J, Chi X, Yang Z. Bridging traditional and machine learning-based algorithms for solving PDEs: the random feature method. J Mach Learn, 2022, 1: 268-98.

\bibitem{CW2003}Chen, Z., Wu, H.: An adaptive finite element method with perfectly matched absorbing layers for the wave scattering by periodic structures. SIAM J. Numer. Anal. 41, 799–826 (2003)

\bibitem{PWNN_2022}Cui T, Wang Z, Xiang X. An efficient neural network method with plane wave activation functions for solving Helmholtz equation. Computers \& Mathematics with Applications, 2022, 111: 34-49.

\bibitem{RobustTrain_2020}Cyr E C, Gulian M A, Patel R G, et al. Robust training and initialization of deep neural networks: An adaptive basis viewpoint. Mathematical and Scientific Machine Learning. PMLR, 2020: 512-536.

\bibitem{DRM_2018} E W, Y B. The deep Ritz method: A deep learning-based numerical algorithm for solving variational problems. Commun Math Stat, 2018, 6(1): 1-12.

\bibitem{ABC_1977}Engquist B, Majda A. Absorbing boundary conditions for numerical simulation of waves. Proceedings of the National Academy of Sciences, 1977, 74(5): 1765-1766.

\bibitem{Approximator1}Kurt Hornik, Maxwell Stinchcombe, Halbert White, Multilayer feedforward networks are universal approximators, Neural Netw. 2 (5) (1989) 359–366.

\bibitem{Approximator2} Kurt Hornik, Maxwell Stinchcombe, Halbert White, Universal approximation of an unknown mapping and its derivatives using multilayer feedforward networks,Neural Netw. 3 (5) (1990) 551–560.

\bibitem{TBC_2024}Lin L, Lv J, Li S. An adaptive finite element DtN method for the acoustic-elastic interaction problem. advances in computational mathematics, 2024, 50(4): 67.

\bibitem{TBC_2025}Lin L, Lv J, Niu T. An adaptive DtN-FEM for the scattering problem from orthotropic media. Applied Numerical Mathematics, 2025, 209: 140-154.

\bibitem{PINN_2019} Raissi M, Perdikaris P, Karniadakis G E. Physics-informed neural networks: A deep learning framework for solving forward and inverse problems involving nonlinear partial differential equations. Journal of Computational physics, 2019, 378: 686-707.

\bibitem{JLLWWZ2022}Xue Jiang, Peijun Li, Junliang Lv, Zhoufeng Wang, Haijun Wu and Weiying Zheng, An adaptive edge finite element DtN method for Maxwell’s equations in biperiodic structures, IMA Journal of Numerical Analysis, 42(3), 2794-2828, 2022.

\bibitem{JLLZ2017_ESAIM}Xue Jiang, Peijun Li, Junliang Lv and Weiying Zheng, An adaptive finite element PML method for the elastic wave scattering problem in periodic structures, ESAIM: Mathematical Modelling and Numerical Analysis, 51(5), 2017-2047, 2017.

\bibitem{JLLZ2017_JSC}Xue Jiang, Peijun Li, Junliang Lv and Weiying Zheng, An adaptive finite element method for the wave scattering with transparent boundary condition, Journal of Scientific Computing, 72(3), 936-956, 2017.

\bibitem{TBC_2021}Xu L, Yin T. Analysis of the Fourier series Dirichlet-to-Neumann boundary condition of the Helmholtz equation and its application to finite element methods. Numerische Mathematik, 2021, 147(4): 967-996.

\bibitem{PAL_2019}Yang Z, Wang L L, Gao Y. A Truly Exact and Optimal Perfect Absorbing Layer for Time-harmonic Acoustic Wave Scattering Problems. arXiv preprint arXiv:1910.08884, 2019.

\bibitem{RBDNN_2023}Yang A L. A novel deep neural network algorithm for the Helmholtz scattered problem in the unbounded domain. International Journal of Numerical Analysis \& Modeling, 2023, 20(5).

\bibitem{RBDNN_2022}Yang A L, Gu F. Mesh-less, ray-based deep neural network method for the Helmholtz equation with high frequency. International Journal of Numerical Analysis \& Modeling, 2022, 19(4).

\bibitem{WAN_2020}Zang Y, Bao G, Ye X, et al. Weak adversarial networks for high-dimensional partial differential equations. Journal of Computational Physics, 2020, 411: 109409.

\bibitem{cPINN_2024}Zhang R, Gao Y. Learning scattering waves via coupling physics-informed neural networks and their convergence analysis. Journal of Computational and Applied Mathematics, 2024, 446: 115874.

\bibitem{TransNet_2024}Zhang Z, Bao F, Ju L, et al. Transferable Neural Networks for Partial Differential Equations. Journal of Scientific Computing, 2024, 99(1): 2.

\bibitem{SNN_arXiv}Zhaodong Xu, Zhiqiang Sheng. Subspace method based on neural networks for solving the partial differential equation. arXiv:2404.08223.

\end{thebibliography}
\end{document}